\documentclass[fleqn,10pt]{wlscirep}
\usepackage[utf8]{inputenc}
\usepackage[T1]{fontenc}
\usepackage{subcaption}
\usepackage{amsthm}
\newtheorem{thm}{Theorem}
\newtheorem{definition}{Definition}
\newtheorem{prop}[thm]{Proposition}
\newtheorem{lem}[thm]{Lemma}
\title{On some properties of the compliance-volume fraction Pareto front in topology optimization useful for material selection.}

\author[1,*]{Edouard Duriez}
\author[1]{Miguel Charlotte}
\author[2]{Catherine Azzaro-Pantel}
\author[1]{Joseph Morlier}
\affil[1]{ICA, Université de Toulouse, ISAE-SUPAERO, MINES ALBI, UPS, INSA, CNRS,3 Rue Caroline Aigle, 31400 Toulouse (France)}
\affil[2]{Laboratoire de Génie Chimique, Université Toulouse, CNRS, INPT, UPS, Toulouse (France)}

\affil[*]{edouard.duriez@gmail.com}

\begin{abstract}
Selecting the optimal material for a part designed through topology optimization is a complex problem. The shape and properties of the Pareto front plays an important role in this selection. In this paper we show that the compliance-volume fraction Pareto fronts of some topology optimization problems in linear elasticity share some useful properties. These properties provide an interesting point of view on the efficiency of topology optimization compared to other design approaches such as parametric structural optimization. We construct a simple meta-model which requires only one full topology optimization to fit the whole Pareto fronts. Precise Pareto fronts are obtained independently. The fast meta-model constructed has a maximum error of 6.4\% with respect to these precise Pareto fronts, on the different problems tested. The selection of the optimal material is then successfully tested on the mass minimization of an MBB beam with an illustrative choice of 4 materials.
\end{abstract}
\begin{document}

\flushbottom
\maketitle
%
%
\thispagestyle{empty}


\section*{Introduction}

Classical topology optimization consists in optimizing a structure performance by distributing material within a design space with a free choice of topology (\cite{bendsoe_topology_2004}). There is a great variety of methods in this research field: homogenization (\cite{bendsoe_generating_1988}), Solid Isotropic Material with Penalization (SIMP) (\cite{bendsoe_bendsoe_1989}), evolutionary methods (\cite{xie_simple_1993}), level set methods (\cite{wang_level_2003}), moving morphable components (MMC) (\cite{guo_doing_2014}), generalized geometry projection(GGP) (\cite{coniglio_generalized_2019}), among others. These methods are described in more detail in \cite{xia_bi-directional_2018} and \cite{norato_topology_2018}.

Let us consider a classical topology optimization framework, formulated as a material distribution problem such as the one in Eq. (\ref{topoptpb}). The design domain $\Delta$ is discretized into elements $e$ characterized by their density $\rho_e$. The compliance $C$ of the design, defined by the node displacements $\textbf{U}$ and the stiffness matrix $\textbf{K}$, is minimized under a constraint on volume fraction $V_f$. The design is subject to a load $\textbf{F}$.

\begin{subequations}
\begin{align}
    \underset{\rho_e}{\arg{}\min{}} C(\rho_e) &= \textbf{U}^T \textbf{K} \textbf{U} \\
    s.t.\quad & \textbf{F}=\textbf{K} \textbf{U} \\
    & 0  \leq \rho_e \leq 1, \quad \forall e \in \Delta \\
    & \underset{e}{\sum} \rho_e \leq V_f
\end{align}
\label{topoptpb}
\end{subequations}

\begin{figure}
  \centering
  \includegraphics[width=0.5\textwidth]{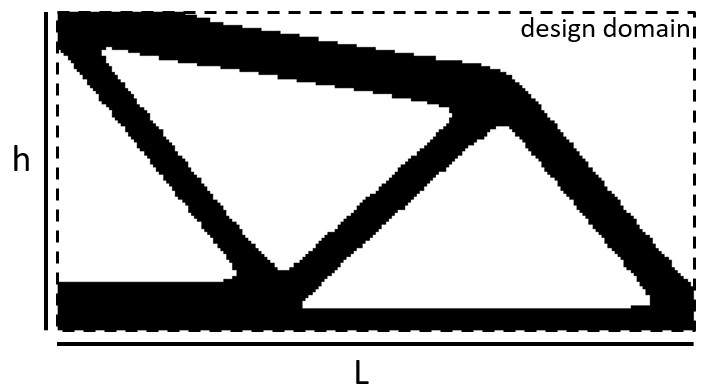}
\caption{Example of topology obtained at the end of an optimization problem associated to Eq. (\ref{topoptpb})}
\label{topol}       
\end{figure}

In this paper, we use the SIMP approach through the well-known code \emph{top88.m}(\cite{andreassen_efficient_2011}) to solve this problem. In this approach, during the optimization process, each element $e\in\Delta$ has an assigned Young modulus $E_e$ obeying Eq. (\ref{penal}). We use a penalisation value $p$ of 3.
\begin{equation}
    E_e = \rho_e^p E_{r}
\label{penal}
\end{equation}
Here $E_{r}$ stands as a reference Young modulus. For illustration, an example of design or topology obtained after the optimisation can be seen in Fig. \ref{topol}.

\begin{figure}
  \centering
  \includegraphics[width=0.7\textwidth]{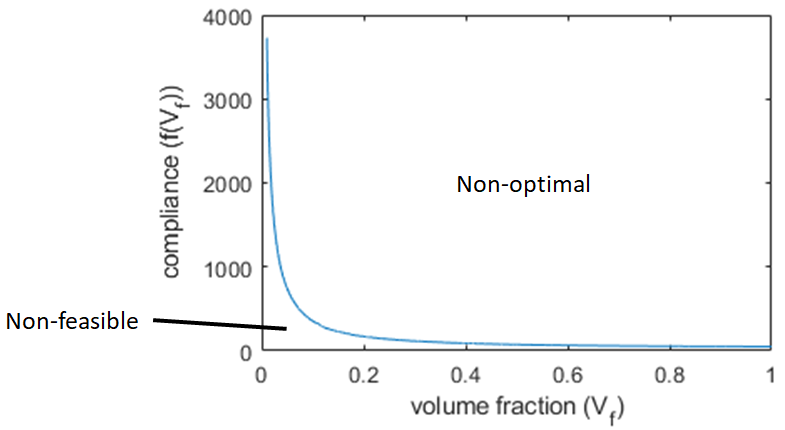}
\caption{MBB beam Pareto front: optimal compliance ($C_{opt}$) as a function of volume fraction.}
\label{Paretop1}       
\end{figure}

Besides, we can carry out an optimization for each discrete value of volume fraction, and plot the final optimal compliance $C_{opt}$ as a function of the volume fraction $V_f$, as in Fig. \ref{Paretop1}. This plot is so the Pareto front of the \emph{multi-objective topology optimization} problem for $C$ and $V_f$.

\cite{fu_recent_2020} reviews multi-objective topology optimization. In order to speed up the computation of the compliance-volume fraction Pareto front, it is possible to start with a high volume fraction and travel along the Pareto front from one design to the next one with decreasing volume fraction. This can be achieved using topological derivatives (\cite{turevsky_efficient_2011}, \cite{suresh_199-line_2010}). \cite{novotny_topological-shape_2006} provides useful expressions for the topological-Shape sensitivity. This method is very fast, but is sometimes far from the global optimum, in particular for low densities, as shown in Section \ref{numPareto}. 

In this paper we will use compliance-volume fraction Pareto fronts to solve the problem of the mass $M$ minimization of a part with a design $\mathcal{D}$ and material $mat$ as variables. The optimization is subject to a constraint $\delta_{max}$ imposed on the deflection $U_F$ generated by a prescribed loading force $F$ on its point of application. The materials are defined by their Young's modulus ($E$) and their density $\rho$. They can be chosen among a set of materials $\Phi$. The part can have any design in the design domain. The volume fraction of the design, or ratio between the volume of the part and the volume of the design domain, $V_f(\mathcal{D})$ therefore stays in [0,1]. This problem is summarized in Eq. (\ref{galpb}).

\begin{subequations}
\begin{align}
    \underset{mat,\mathcal{D}}{\arg{}\min{}} \quad & M (mat,\mathcal{D}) \\
    s.t.\quad & U_F \leq \delta_{max} \\
    & mat = \{E, \rho \} \in \Phi \\
    & 0 < V_f(\mathcal{D}) \leq 1
\end{align}
\label{galpb}
\end{subequations}

It should be noted that Poisson's ratio is not included as a material variable. This choice is made for the sake of simplicity, and since there is only a limited influence on the compliance obtained via topology optimization over the range of Poisson's ratio for typical materials. In order to take Poisson's ratio into account, more complex meta-models than the one built in Section \ref{buildmeta} would have to be created. Only 2D topology optimization examples are shown in this article, but the results can easily be applied for 3D.

Now, since the material used and the structure are coupled, they must be optimized simultaneously through multi-disciplinary design optimization (MDO) (\cite{martins_multidisciplinary_2013},\cite{martins_engineering_2021}) and this can be computationally intensive. An alternative is to project the discrete materials on a continuous space and proceed with a neural network (\cite{suresh_integrating_2021}). For simpler problems, it Ashby indices can be used (\cite{ashby_materials_2004}), which are described and used hereafter. 

Ashby's method aims at simplifying material selection. It can help uncouple design optimization from material selection. The general steps of the method are as follows. First, the objective $P$ is written as a function $f_1$ of free geometrical variables $G$, fixed parameters $F$, and material variables $M$:
\begin{equation}
    P=f_1(G,F,M)\,.
\end{equation}
As a second step, the free geometrical variables $G$ are replaced by fixed parameters and material variables, by using the constraints:
\begin{equation}\label{F1F2Eq}
    P=f_1(f_2(F,M),F,M) \,.
\end{equation}
If the new expression in Eq.~\ref{F1F2Eq} of the objective function  is separable into a product of a function $f_3$ of the fixed parameters $F$ by a function $f_4$ of the material variables $M$, then $f_4$ is named the material index:
\begin{equation}
    P=f_3(F)\times  f_4(M)\,.
\end{equation}
Indeed, selecting the material that maximizes (or minimizes) this index is equivalent to maximizing (or minimizing) the objective function and this independently of the values of the other terms (\cite{Ashby2002MaterialsAD}).

Classical Ashby indices are used for pre-defined problems such as: generic components with simple section shapes such as a tie; a tensile component, a panel loaded in bending or a beam loaded in bending. These indices have been extended to the design of trusses (\cite{ananthasuresh_concurrent_2003},\cite{rakshit_simultaneous_2007}). To our knowledge, these indices have only been used once with topology optimization, in the specific case of compliant mechanism optimization (\cite{achleitner_material_2022}).

An objective of this work is to rethink Ashby's indices in order to minimize the mass of structures optimized through topology optimization under a stiffness constraint (equivalent to a displacement constraint if the load is imposed), as in Eq. (\ref{galpb}). The case where the out-of-plane thickness of the design domain is variable has been addressed in \cite{duriez_ecodesign_2022}. In that case, the objective function (mass) is separable in an expression independent of the material multiplied by the ratio $\frac{\rho}{E}$. The optimal material is therefore the one with the lowest $\frac{\rho}{E}$ ratio. If the thickness $t$ is fixed however, the problem is much more complex, as the objective function is no longer separable. The mass can be written as in Eq. (\ref{mass}), where $L$ and $h$ are the lengths and heights of the design domain, and $V_f$ is the volume fraction of the design.
\begin{equation}
    M=L h t V_f \rho
\label{mass}
\end{equation}
The compliance of the structure is $C=\displaystyle{\frac{f(V_f) F^2}{t E}}$, with $E$ the material's Young Modulus and $f(V_f)$ the compliance of the optimized design obtained for a unit force, a unit thickness, a unit Young's modulus and a volume fraction $V_f$. This function $f$ is the Pareto front highlighted in Fig. \ref{Paretop1}. Let us consider for simplicity's sake that the load is applied in only one node and has a value of $F$. In this case, we have $C=\textbf{U}^T\textbf{F}\equiv U_F F$.

For the optimal design, the stiffness constraint is activated and consecutively the constraint can be written as in Eq. (\ref{constr})
\begin{equation}
    \delta_{max}=U_F = \frac{C}{F} = \frac{f(V_f) F}{t E}
\label{constr}
\end{equation}
When $f$ is a continuous and (as proved in Section \ref{npos}) monotonically decreasing function of $V_f$, we can use so its reciprocal $f^{-1}$ to get an expression of the volume fraction $V_f$:
\begin{equation}
    V_f=f^{-1}(\frac{t E \delta_{max}}{F})
\label{geomvar}
\end{equation}
This in turn enables us to remove the geometric free variable $V_f$ in the expression of the objective function:
\begin{equation}
    M=L  h  t  f^{-1}(\frac{t E \delta_{max}}{F}) \rho
\label{mass2}
\end{equation}
This expression is not separable in general, but \begin{equation}
f_4 \stackrel{\rm def}{=} f^{-1}(\frac{t E \delta_{max}}{F}) \rho    \label{AshbyIndex}
\end{equation} 
can be considered as our Ashby index.

In order to be able to select the optimal material, this paper is organized as follows. The first section \ref{proof} presents some properties of the compliance-volume fraction Pareto front  through a new function called ER (Efficiency Ratio). We then develop an original meta-model to obtain fast and reliable Pareto fronts in a second section. Finally, we introduce the concurrent method to select the optimal material and design.

\section*{Proof that the ER always lies in the interval [0,1]}
\label{proof}
\subsection*{Some general statements on topology optimization}
The aim of this proof is to be able to limit the choice of materials for a structure obtained through topology optimization to a small set.
In the remainder of this proof we will notably consider a theoretical case where the domain is discretized into an infinite number of infinitely small elements. This means that infinitely small features can appear. 

We consider theoretically that a topology optimization is carried out for every volume fraction in $]0,1]$. We consider that the global optimum is reached for each of these optimizations. If there is no global optimum, the tight upper bound can be considered instead (\cite{allaire_chapter_2021}). 

\begin{figure}
  \centering
  \includegraphics[width=0.7\textwidth]{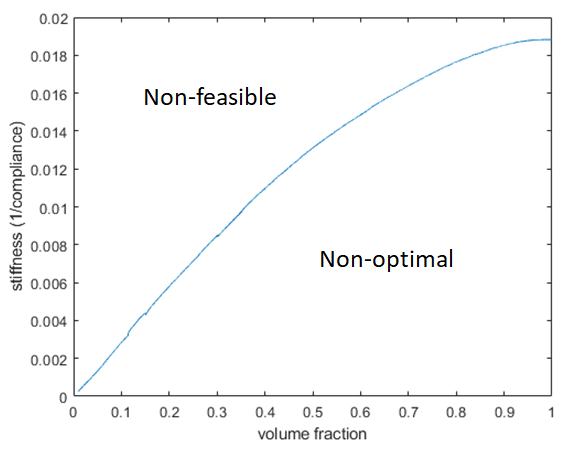}
\caption{MBB beam Pareto front: optimal stiffness ($\kappa_{opt}$) as a function of volume fraction, for the same data as in Fig. \ref{Paretop1}.}
\label{Paretostiff}       
\end{figure}

By \emph{optimal stiffness} ($\kappa_{opt}$) we mean the inverse of the optimal compliance, such that $\kappa_{opt}(V_f)=1/C_{opt}(V_f)$.
We can plot these optimal stiffnesses  as a function of the volume fraction, as in Fig. \ref{Paretostiff}. This plot is another Pareto-front of the same problem as in Fig. \ref{Paretop1}. The optimal designs can change very suddenly and completely along the Pareto front, as illustrated in Fig. \ref{angular}, but the volume fraction is continuous.

Let us consider the derivative $\displaystyle{\frac{d \kappa_{opt}}{d V_f}}(V_f)$ of this function $\kappa_{opt}(V_f)$. This derivative gives us the increase in stiffness due to a marginal increase in volume. Therefore, for every  value $V_f$, this derivative corresponds to the contribution of the last volume that was added to the design, starting from a void design, and going up the Pareto front. We name this derivative the efficiency of volume addition. The mean value of the efficiency of all the material added to attain a given design can consecutively be introduced like
\begin{equation}
    \overline{\frac{d \kappa_{opt}}{d V_f}} \stackrel{\rm def}{=} \frac{1}{V_f-0} \int_{0}^{V_f} \frac{d \kappa_{opt}}{d V_f}(V_f) dV_f =\frac{\kappa_{opt}}{V_f} 
\label{meaneff}
\end{equation}
at the condition that $\kappa_{opt}$ is a piecewise continuously differentiable function. We suppose this condition to be true. This is verified numerically in Section \ref{numPareto}.

We define the ``efficiency ratio'' (ER) of a design (noted $n(V_f)$) as the ratio between the present efficiency of volume addition and its mean value for all the volume in the design. We obtain for $n$ the expression in Eq. (\ref{n1}).
\begin{equation}
    n(V_f) \stackrel{\rm def}{=} \frac{\displaystyle{\frac{d \kappa_{opt}}{d V_f}}}{\displaystyle{\frac{\kappa_{opt}}{V_f}}} = \frac{\displaystyle{\frac{-C_{opt}'(V_f)}{C_{opt}^2(V_f)}}}{\displaystyle{\frac{1}{V_f C_{opt}(V_f)}}}
    = -V_f \frac{C_{opt}'(V_f)}{C_{opt}(V_f)} 
\label{n1}
\end{equation}
This index $n$ can be interpreted as a measure of how much material added to the design and leading to an increased volume fraction, will be efficient in lowering compliance, compared to the mean efficiency of the material already used. Thus, notably if $n=1$, the material added in small quantity has the same efficiency as the material already used in the design; if $n>1$, the added material has a higher efficiency than the material already in the design.

Let us consider a special case. If the ER is a constant throughout the considered volume fraction range, Eq. (\ref{n1}) can be rearranged and Eqs. (\ref{nconst}) are verified successively, with a nonzero real constant $A>0$.

\begin{subequations}
\begin{align}
    C_{opt}'(V_f) = -n \frac{C_{opt}(V_f)}{V_f} \\
    C_{opt}(V_f) = A V_f^{-n} \mbox{ so that }   \kappa_{opt}(V_f) = \frac{1}{A} V_f^{n}\,.
\end{align}
\label{nconst}
\end{subequations}

\begin{figure}
\begin{subfigure}{.33\textwidth}
  \includegraphics[width=.99\linewidth]{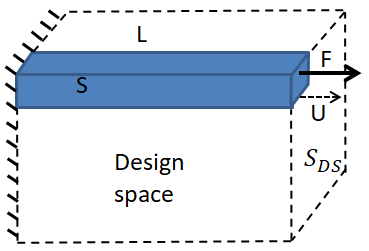}
  \caption{Tension rod: n=1}
\end{subfigure}%
\begin{subfigure}{.33\textwidth}
  \includegraphics[width=.99\linewidth]{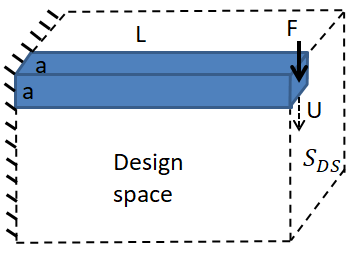}
  \caption{Bending beam: n=2}
\end{subfigure}%
\begin{subfigure}{.33\textwidth}
  \includegraphics[width=.99\linewidth]{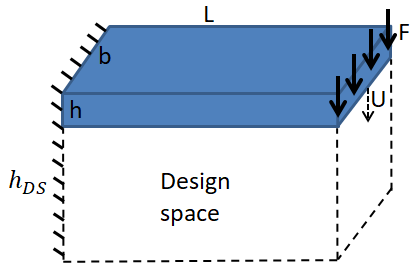}
  \caption{Bending plate: n=3}
\end{subfigure}%
\caption{Simple mechanical components with a constant n.}
\label{simple}
\end{figure}

A constant ER can be found in very simple mechanical components such as a tension rod  ($n=1$) (Eq. (\ref{rod})), a Euler-Bernoulli cantilever bending beam ($n=2$) with a square cross-section (Eq. (\ref{rod})), or a Love-Kirchhoff bending plate ($n=3$) (Eq. (\ref{rod})).
Indeed, in these cases, the apparent stiffnesses have the following expression:
\begin{subequations}
\begin{align}
    \kappa_{rod}(V_f) = \frac{F}{U} = \frac{E S}{L} = \frac{E V_f S_{DS}}{L}; \quad V_f=\frac{S}{S_{DS}} \label{rod}\\
    \kappa_{beam}(V_f) = \frac{F}{U} = \frac{E a^4}{4 L^3} = \frac{E V_f^2 S_{DS}^2}{4 L^3}; \quad V_f=\frac{a^2}{S_{DS}}  \label{beam}\\ 
    \kappa_{plate}(V_f) = \frac{F}{U} = \frac{E b h^3}{4 L^3} = \frac{E b V_f^3 h_{DS}^3}{4 L^3}; \quad V_f=\frac{h}{h_{DS}} \label{plate}
\end{align}
\label{stiffnesses}
\end{subequations}
with $F$ the tip load, $U$ the tip displacement, $E$ the material's Young's modulus, $L$ the length of the rod/beam/plate, $S_{DS}$ the cross-section of the design domain, $S$ the rod's cross-section, $a$ the square beam's thickness, $h$ the plate's thickness, $b$ the plate's width, and $h_{DS}$ the thickness of the design domain volume.
These cases appear in Fig. \ref{simple}. Similar problems are considered in \cite{ashby_materials_2004}.

In the general case of structures optimized through topology optimization, $n$ is not necessarily a constant and  we aim to prove in this section that for such structures, Eq. (\ref{limn}) is always verified
\begin{equation}
    0 \leq n \leq 1\,.
\label{limn}
\end{equation}
Our proof takes the form of three steps which are developed in the three following subsections 
with the following outline: step 1 shows that $\lim\limits_{V_f \searrow 0} n(V_f) = 1$, step 2 shows that we always have $0 \leq n$, and step three proves that $n$ is a monotonous decreasing function of $V_f$.

\subsection*{Proof that the limit of $n(V_f)$ in 0 is 1.}
At very low volume fractions, a good possible design is a truss (\cite{sigmund_non-optimality_2016}) with very thin members joining the load application points and the fixed boundary conditions (if the problem is well-posed: fixed boundary conditions prevent rigid-body motion). The optimal truss design will have its members keep proportional thicknesses to one another, as the total volume fraction varies. As a result, as the volume fraction $V_f$ tends towards 0, the stiffness of the total structure will diminish in proportion with the members' thicknesses, themselves proportional to the volume fraction. In this case the optimal stiffness $\kappa_{truss}$ can be written as in Eq. (\ref{stif0}) with a specific nonzero real constant $\alpha$.
\begin{equation}
    \kappa_{truss}=\alpha \, V_f \quad , \quad \frac{d \kappa_{truss}}{d V_f} = \alpha \quad , \quad \alpha > 0\quad \mbox{ for }V_f>0\,.
\label{stif0}
\end{equation}
If this truss is not the optimal design, the optimal design must be stiffer. \\
As all the designs for $V_f=0^+$ are the same with $\kappa(0^+)=0$ ($C(0^+)=\infty$), we have,
\begin{equation}
    \frac{d \kappa_{opt}}{d V_f}(0^+) \geq \frac{d \kappa_{truss}}{d V_f}(0^+) = \alpha  > 0\,.
\label{derivstif0}
\end{equation}

Because $\kappa_{opt}(0^+) = 0$, we can write Eq. (\ref{deriv}).
\begin{equation}
    \lim_{V_f \searrow 0} \frac{\kappa_{opt}(V_f)}{V_f} = \lim\limits_{V_f \searrow 0} \frac{\kappa_{opt}(V_f) - \kappa_{opt}(0^+)}{V_f - 0} = \frac{d \kappa_{opt}}{d V_f}(0^+) > 0
\label{deriv}
\end{equation}

Therefore, we can write Eq. (\ref{limn1}). 
\begin{equation}
    \lim_{V_f \searrow 0} n = \frac{\lim\limits_{V_f \searrow 0} \frac{d \kappa_{opt}}{d V_f}(V_f)}{\lim\limits_{V_f \searrow 0} \frac{\kappa_{opt}(V_f)}{V_f}} = \frac{\frac{d \kappa_{opt}}{d V_f}(0^+)}{\frac{d \kappa_{opt}}{d V_f}(0^+)} = 1
    \, .
\label{limn1}
\end{equation}

\subsection*{Proof that n is positive.}
\label{npos}
Let us prove that $\kappa_{opt}$ is an increasing function. We proceed by contradiction.
\\
We suppose that $\kappa_{opt}$ is not an increasing function of $V_f$.
Therefore,
\begin{center}
$\exists V_f \in ]0,1]$, $\exists dV_f>0$ such that $\kappa_{opt}(V_f) > \kappa_{opt}(V_f + dV_f)$.
\end{center}
We consider the design obtained at the end of the topology optimization for this volume fraction $V_f$, and the stiffness of which is $\kappa_{opt}(V_f)$. We can increase the density of all elements with $\rho \neq 1$ by a same amount such that we obtain a new design with a volume fraction increased by $dV_f$. This new design is stiffer than the initial one, because of Eq. (\ref{penal}). Therefore, we have found a design with a volume fraction of $V_f+dV_f$ and a higher (effective) stiffness ($\kappa_1$) than the initial (optimal) one at $V_f$. This translates into Eq. (\ref{compare}).
\begin{equation}
    \kappa_1 > \kappa_{opt}(V_f) > \kappa_{opt}(V_f + dV_f)\,.
\label{compare}
\end{equation}
However, $\kappa_{opt}(V_f + dV_f)$ is the highest stiffness accessible to a design of volume fraction $V_f+dV_f$, because it is on the Pareto front. 
Therefore, \mbox{$\kappa_1 \leq \kappa_{opt}(V_f + dV_f)$}. 
This is an evident contradiction with Eq. (\ref{compare})
and entails therefore that $\kappa_{opt}$ is an increasing function of $V_f$ and $\forall V_f \in ]0,1], \kappa_{opt}'(V_f) \geq 0$.

Furthermore, we know that $0 < V_f \leq 1$ and $\kappa_{opt}(V_f) >0$.  
Hence, as a conclusion,
\begin{equation}
    n(V_f) = \frac{\displaystyle{\frac{d \kappa_{opt}}{d V_f}}}{\displaystyle{\frac{\kappa_{opt}}{V_f}}} \geq 0\,,\quad \forall V_f \in ]0,1]
\label{conc1}
\end{equation}

\subsection*{Proof that $n \leq 1$} 
In this subsection we first prove the above assertion in a particular case where the topologies vary sufficiently continuously. Then, we explain why we assume whether this can be generalized or not.
\subsubsection*{Any angular point on $\kappa_{opt}$ due to a small change of topology can only be concave.}
\begin{figure}
\begin{subfigure}{.5\textwidth}
  \includegraphics[width=.99\linewidth]{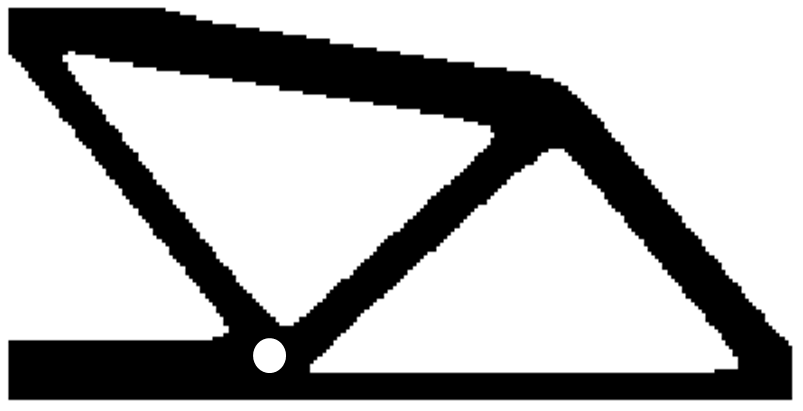}
  \caption{The  removed (white) disc has a volume of $\delta$.}
\end{subfigure}%
\begin{subfigure}{.5\textwidth}
  \includegraphics[width=.99\linewidth]{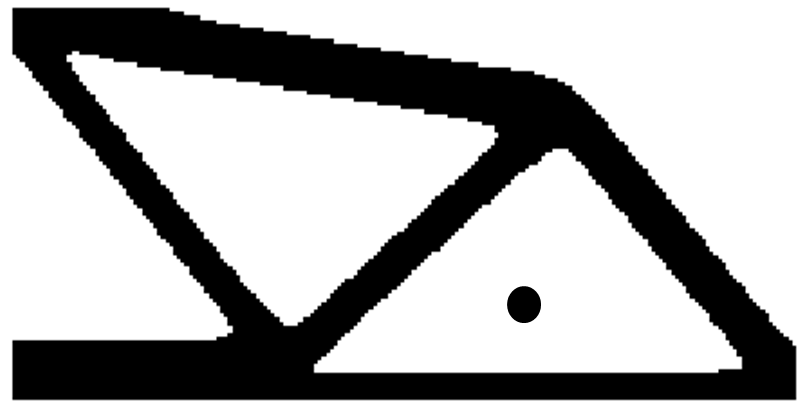}
  \caption{The added (black) disc has a volume of $\delta$.}
\end{subfigure}%
\newline
\begin{subfigure}{.5\textwidth}
  \includegraphics[width=.99\linewidth]{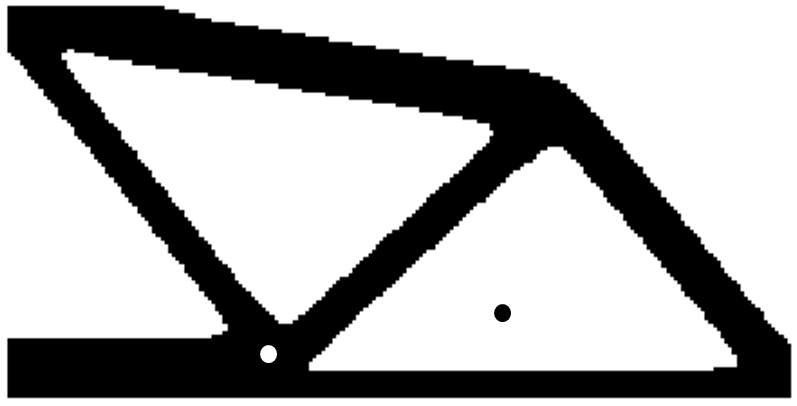}
  \caption{The added and removed discs both have a volume of $\delta/2$.}
\end{subfigure}%
\begin{subfigure}{.5\textwidth}
  \includegraphics[width=.99\linewidth]{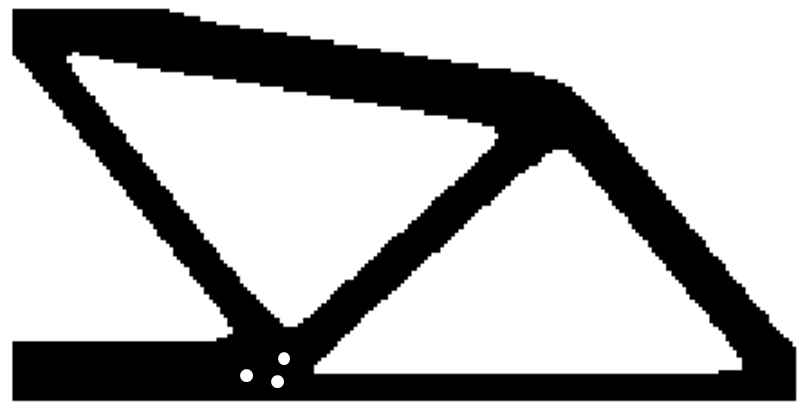}
  \caption{The 3 removed discs all have a volume of $\delta/3$.}
\end{subfigure}%
\caption{Different topologies that are $\delta$-appart from the topology in Fig. \ref{topol}.}
\label{deltaappart}
\end{figure}

We use the definitions and the Lemma provided in \cite{turevsky_efficient_2011} to work with topologies. We recall them below:

\begin{definition}
Topologies $\Omega$ and $\Omega'$ are $\delta$-apart if their symmetric volume difference is less than $\delta$, i.e., 
\begin{equation}
\Delta V(\Omega,\Omega')=|\Omega \setminus \Omega'|+|\Omega \setminus \Omega'| \leq \delta
\label{apart}
\end{equation}
\end{definition}
Fig. \ref{deltaappart} illustrates various topologies that are $\delta$-apart from Fig. \ref{topol}.

\begin{definition}
A topology $\Omega$ is said to be locally Pareto-optimal, if it is Pareto-optimal with respect to all topologies that are within a distance $\delta$ apart from it, where $\delta$ is sufficiently small, i.e.
\begin{equation}
\kappa(\Omega) \geq \kappa(\Omega') \quad \forall \Omega',  \Delta V(\Omega,\Omega') \leq \delta, |\Omega'| \leq |\Omega|
\label{localopt}
\end{equation}

\end{definition}

We add our own definition:
\begin{definition}
The set $\mathcal{N}(\Omega,\delta)$ of $\delta$-neighbours of topology $\Omega$ is defined as
\begin{equation}
\mathcal{N}(\Omega,\delta) \stackrel{\rm def}{=} \{ \Omega', \Delta(\Omega,\Omega') \leq \delta \}
\label{neighbours}
\end{equation}
\end{definition}



The definitions of topological sensitivity from \cite{novotny_topological-shape_2006} are recalled below:

\begin{definition}
The subtraction topological sensitivity $\mathcal{T}^S$ is defined as
\begin{equation}
\mathcal{T}^S(p) \stackrel{\rm def}{=} \lim_{\epsilon \searrow 0} \frac{F(\Omega \setminus B_{\epsilon}(p))-F(\Omega)}{|B_{\epsilon}(p)|}
\label{subsensitivity}
\end{equation}
where $B_{\epsilon}(p)$ is a ball of radius $\epsilon>0$ and center p$\in \Omega$ and $F$ is a cost function to be minimized under a constraint on the volume $|\Omega|$.
\end{definition}

$\mathcal{T}^S$ therefore represents the change in cost function value when an infinitely small hole is drilled in the design. It is therefore always positive. Note that $\mathcal{T}^S$ is a field defined in each point on the topology.
In the same way, the additive topological sensitivity $\mathcal{T}^A$ is defined as
\begin{equation}
\mathcal{T}^A(q) \stackrel{\rm def}{=} \lim_{\epsilon \searrow 0} \frac{F(\Omega \cup B_{\epsilon}(q))-F(\Omega)}{|B_{\epsilon}(q)|}
\label{addsensitivity}
\end{equation}
with the same notation and $q \in D \setminus \Omega$, where D is the design space. It therefore represents the change in cost function when an infinitely small point of matter is added outside the design. It is therefore always negative.

\begin{figure}
  \centering
  \includegraphics[width=0.7\textwidth]{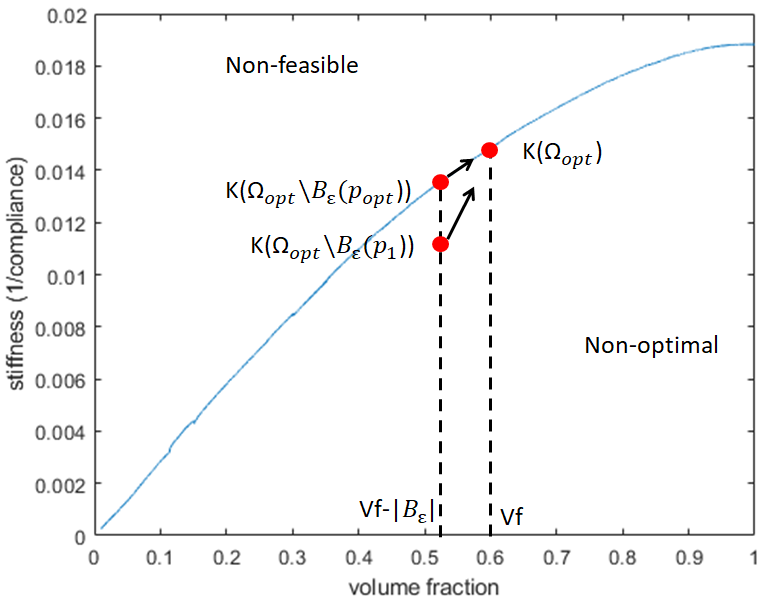}
\caption{The stiffness is plotted as a function of the volume fraction along an example of local Pareto front. The stiffness of $\Omega$ is plotted, as well as $\Omega \setminus B_{\epsilon}(p)$ for two points: $p_{opt}$ corresponding to the minimum of $\mathcal{T}^S$ and another non-optimal point $p_1$.
The left derivative $\displaystyle{\frac{\partial \kappa}{\partial V_f^-}}$ on the local Pareto front is the minimum $min(\mathcal{T}^S)$ of the subtractive topological derivative.}
\label{derivtop}       
\end{figure}

Let us consider a Pareto-optimal topology $\Omega_{opt}$. We denote $\kappa(\Omega)$ the stiffness of the design associated with topology $\Omega$,  such that $\kappa_{opt}(V_f)\equiv\kappa(\Omega_{opt})$. For a sufficiently small $\delta$, on the local Pareto front on $\mathcal{N}(\Omega_{opt},\delta)$, $\kappa_{opt}(V_f-|B_{\epsilon}|)=\displaystyle{\max_{p\in \Omega_{opt}}}(\kappa(\Omega_{opt} \setminus B_{\epsilon}(p)))$, as can be seen in Fig. \ref{derivtop}. Therefore, the left derivative of $\kappa_{opt}(V_f)$ on the local Pareto front on $\mathcal{N}(\Omega_{opt},\delta)$ is  
\begin{equation}
\renewcommand{\arraystretch}{1.7}\begin{array}{rcl}
\displaystyle{\frac{d \kappa_{opt}}{d V_f^-}} &=& \lim\limits_{\epsilon \searrow 0}  \displaystyle{\frac{\kappa(\Omega_{opt})-\displaystyle{\max_{p\in \Omega_{opt}}}(\kappa(\Omega_{opt} \setminus B_{\epsilon}(p)))}{|B_{\epsilon}(p)|}} \\
&& \equiv \lim\limits_{\epsilon \searrow 0}\ \min\limits_{p\in \Omega_{opt}}\left(\displaystyle{\frac{\kappa(\Omega_{opt})-\kappa(\Omega_{opt} \setminus B_{\epsilon}(p))}{|B_{\epsilon}(p)|}}\right)
\end{array}
\end{equation}
similarly, the right derivative of $\kappa_{opt}(V_f)$ on the local Pareto front on $\mathcal{N}(\Omega_{opt},\delta)$ is
\begin{equation}
\frac{d \kappa_{opt}}{d V_f^+} = \lim_{\epsilon \searrow 0} \max_{q\in D \setminus \Omega_{opt}}(\frac{\kappa(\Omega_{opt} \cup B_{\epsilon}(q))-\kappa(\Omega_{opt})}{|B_{\epsilon}(q)|})
\end{equation}
We suppose appropriate mathematical properties of $\kappa$ and $\kappa_{opt}$, enabling us to write
\begin{subequations}
\begin{align}
\frac{d \kappa_{opt}}{d V_f^-} = \min_{p\in \Omega_{opt}}(\lim_{\epsilon \searrow 0} \frac{\kappa(\Omega_{opt})-\kappa(\Omega_{opt} \setminus B_{\epsilon}(p))}{|B_{\epsilon}(p)|}) \\
\frac{d \kappa_{opt}}{d V_f^+} = \max_{q\in D \setminus \Omega_{opt}}(\lim_{\epsilon \searrow 0} \frac{\kappa(\Omega_{opt} \cup B_{\epsilon}(q))-\kappa(\Omega_{opt})}{|B_{\epsilon}(q)|})
\end{align}
\label{leftrightderiv}
\end{subequations}
as can be seen in Fig. \ref{derivtop}.

\cite{turevsky_efficient_2011} also give us a lemma relating local Pareto-optimality and topological sensitivity:

\begin{lem}
If a topology $\Omega$ is locally Pareto-optimal, then 
\begin{equation}
\min\limits_{q\in D \setminus \Omega}(\mathcal{T}^A(q)) + \min\limits_{p\in\Omega}(\mathcal{T}^S(p)) \geq 0 
\label{Paretoineq}
\end{equation}
\end{lem}

We apply this lemma to $-\kappa$, a cost function to be minimized.
$\Omega_{opt}$ is a Pareto-optimal topology. It is therefore a locally Pareto-optimal topology. Therefore, the following equations are successively obtained:
\begin{subequations}
\begin{align}
\min_{q\in D \setminus \Omega_{opt}}(\lim_{\epsilon \searrow 0} \frac{-\kappa(\Omega_{opt} \cup B_{\epsilon}(q))+\kappa(\Omega_{opt})}{|B_{\epsilon}(q)|}) + \min_{p\in \Omega_{opt}}(\lim_{\epsilon \searrow 0} \frac{-\kappa(\Omega_{opt} \setminus B_{\epsilon}(p))+\kappa(\Omega_{opt})}{|B_{\epsilon}(p)|}) \geq 0
\\
-\max_{q\in D \setminus \Omega_{opt}}(\lim_{\epsilon \to 0} \frac{\kappa(\Omega_{opt} \cup B_{\epsilon}(q))-\kappa(\Omega_{opt})}{|B_{\epsilon}(q)|}) +\min_{p\in \Omega_{opt}}(\lim_{\epsilon \searrow 0} \frac{\kappa(\Omega_{opt})-\kappa(\Omega_{opt} \setminus B_{\epsilon}(p))}{|B_{\epsilon}(p)|}) \geq 0
\\
\frac{d \kappa_{opt}}{d V_f^-} \geq \frac{d \kappa_{opt}}{d V_f^+}
\end{align}
\label{kappaconcave}
\end{subequations}

\begin{figure}
  \centering
  \includegraphics[width=0.7\textwidth]{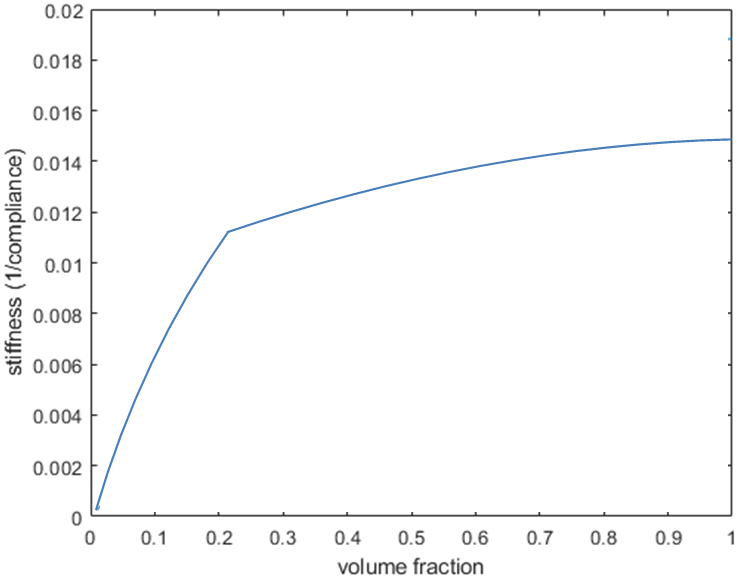}
\caption{Illustrative example of a Pareto front where $\frac{\partial \kappa}{\partial V_f^-} > \frac{\partial \kappa}{\partial V_f^+}$}
\label{concang}       
\end{figure}

Therefore, if a small change in topology creates an angular point in the local Pareto front, it will be a concave angular point as illustrated in Fig. \ref{concang}.

\subsubsection*{Consequence for $n$} 
Provided that $\displaystyle{\frac{d\kappa_{opt}}{dV_f}}$ is a piecewise-continuous, decreasing function of $V_f$ and that \mbox{$\kappa_{opt}(0^+)=0$} as suggested in the foregoing paragraph, 
the following successive equations are then verified for
$\forall x \in ]0,1]$ and $\forall V_f \in ]0,x]$,
\begin{subequations}
\begin{align}
\frac{d\kappa_{opt}}{dV_f}(V_f) &\geq \frac{d\kappa_{opt}}{dV_f}(x) 
\\
\int_0^x \frac{d\kappa_{opt}}{dV_f}(V_f) dV_f \equiv \kappa_{opt}(x) &\geq \int_0^x \frac{d\kappa_{opt}}{dV_f}(x) dV_f \equiv x \frac{d\kappa_{opt}}{dV_f}(x)
\\
1 &\geq  n(x)\equiv\frac{x \frac{d\kappa_{opt}}{dV_f}(x)}{\kappa_{opt}(x)}\,.
\end{align}
\label{intK}
\end{subequations}
Therefore $n(V_f)$ is lower than or equal to 1 for $V_f$ in ]0,1].

\subsubsection*{Remark on the plausible concavity of $\kappa_{opt}$ for 
more general Pareto front types.}
In the previous paragraphs, sufficiently small changes to the topologies along the Pareto front are shown to keep $\kappa$ concave. The proof in \cite{turevsky_efficient_2011} that was used can be adapted to show that 
small shape deformations along the Pareto front also keep $\kappa_{opt}$ concave. Therefore, for a Pareto front 
along which the topologies vary continuously, $\kappa_{opt}$ is concave and $n\leq1$. Pareto fronts built directly by algorithms from one topology to the next  (\cite{suresh_199-line_2010}) fall into this category. 

\begin{figure}
  \centering
  \includegraphics[width=0.7\textwidth]{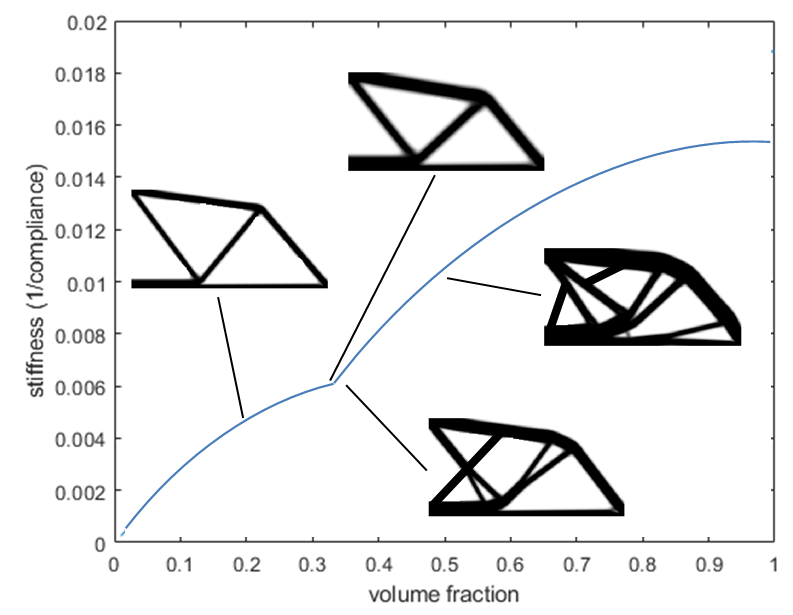}
\caption{Hypothetical illustrative example where a brutal change in topology could result in local convexity.}
\label{angular}       
\end{figure}

However, proving that $\kappa_{opt}$ is concave - if it is true - for more general types of Pareto fronts 
is much harder, and not tackled here. One could indeed imagine for instance that a brutal change in shape or topology could happen in a point on the Pareto front, marking the limit between two concave parts of the Pareto front, but being itself a point of local convexity. This hypothetical case is illustrated in Fig. \ref{angular}. We will however suppose in the remainder of this paper that $n\leq1$ for general compliance-volume fraction Pareto fronts. We observe that $n\leq1$ for all the Pareto fronts we obtain numerically in Section \ref{numPareto}.

\section*{A new meta-model of the Pareto front}
\label{Paretosec}
\subsection*{Numerical strategy to obtain a "precise" Pareto front}
\label{numPareto}
In order to verify numerically the properties of the ER, precise Pareto fronts for different topology optimization problems are needed. These Pareto fronts need to be as close as possible to the real globally optimal Pareto fronts.

\begin{figure}
\begin{subfigure}{.5\textwidth}
  \includegraphics[width=.99\linewidth]{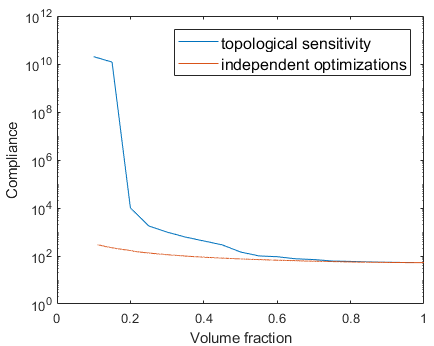}
  \caption{Whole Pareto front.}
\end{subfigure}%
\begin{subfigure}{.5\textwidth}
  \includegraphics[width=.99\linewidth]{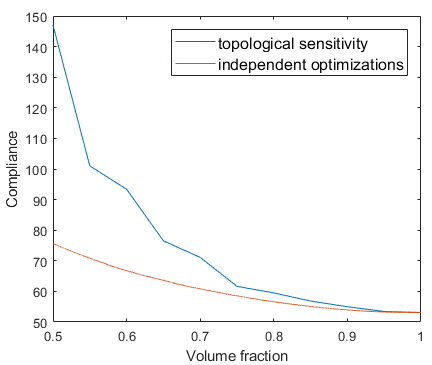}
  \caption{Zoom on high volume fractions.}
\end{subfigure}%
\caption{Comparison between the method from \cite{suresh_199-line_2010}, obtaining the Pareto front from one design to the next through topological sensitivity and a better Pareto front obtained by independent optimizations.}
\label{l199}
\end{figure}

\begin{figure}
  \centering
  \includegraphics[width=0.49\textwidth]{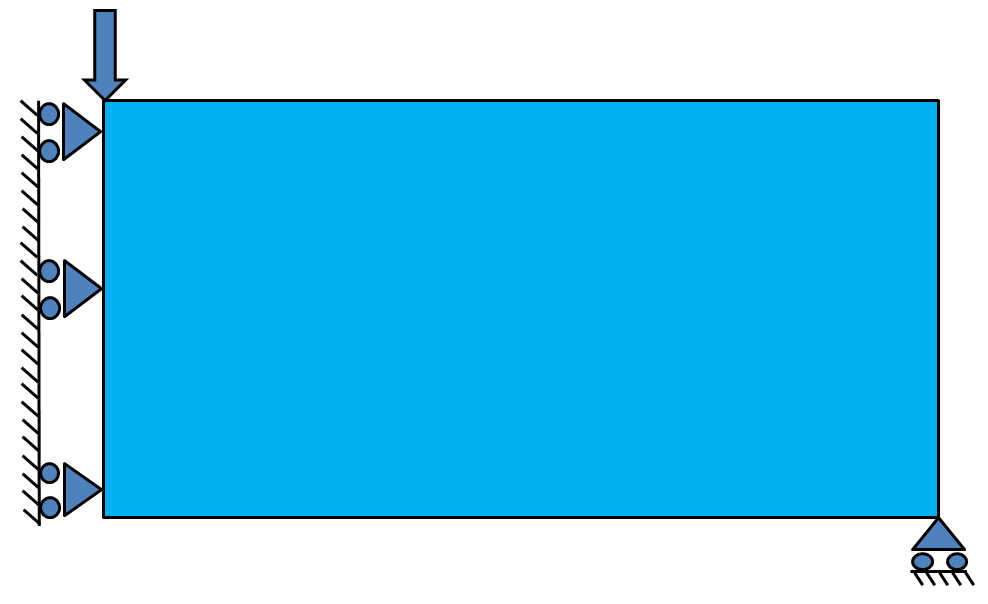}
\caption{The MBB beam problem.}
\label{mbbprob}       
\end{figure}

Obtaining Pareto fronts through typical methods jumping from one design to the next (for example \cite{suresh_199-line_2010}) gives results far from the global optimums. Indeed, big discontinuities in the topologies from one point of the Pareto front to the next are missed with this approach, and it does not work well for low volume fractions. This can be seen in Fig. \ref{l199}, where the method by \cite{suresh_199-line_2010} with volume fraction steps of 0.05, is compared to the method described below on an MBB beam problem (Fig \ref{mbbprob}).

Therefore, we choose to obtain the Pareto front by multiple independent topology optimizations with a varying volume fraction constraint. 10000 regularly spaced points are computed, using the classical SIMP approach implemented in the top88 code (\cite{andreassen_efficient_2011}). We first consider the classical MBB beam problem with a $200 \times 100$ design domain, a penalisation of 3, and a sensitivity/density filter with a radius of 3.
A Pareto front can be obtained by plotting the final compliances of these optimizations. However, this Pareto front would be very noisy. 

\begin{figure}
  \centering
  \includegraphics[width=0.99\textwidth]{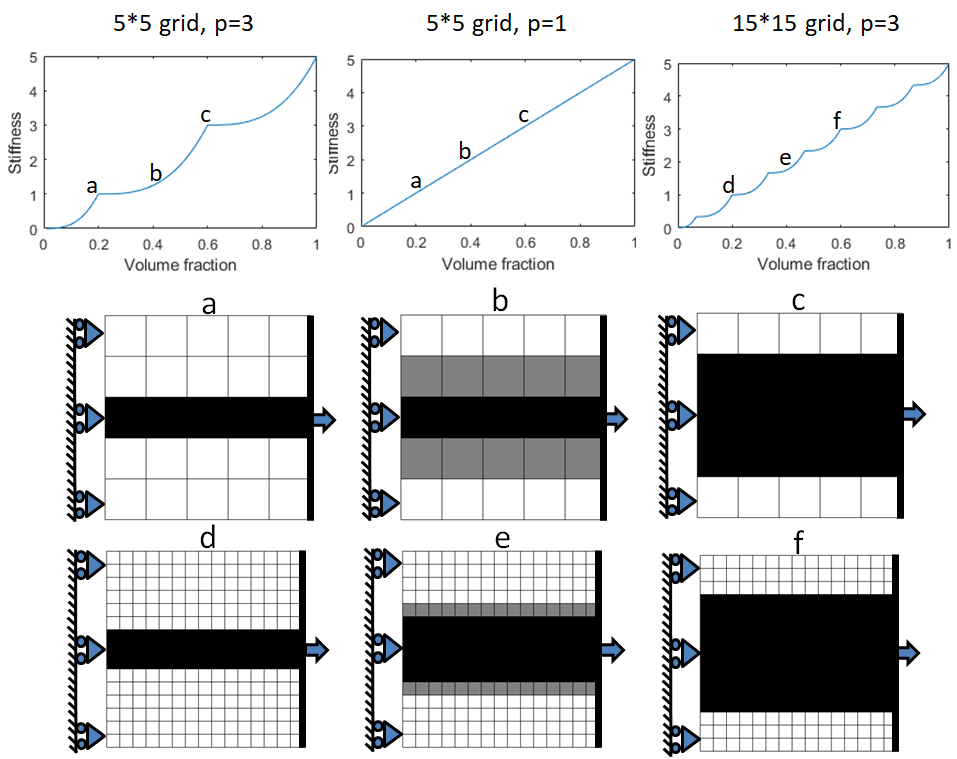}
\caption{Pareto fronts obtained for a simple traction problem, with different penalization (p=1 and p=3) and different discretizations (5*5 or 15*15). The designs for $V_f=0.2$ (a,d), $V_f=0.4$ (b,e) and $V_f=0.2$ (c,f) are given. Using a penalization of 1 to evaluate the design after the optimization gives a Pareto front closer to the one obtained with a higher discretization.}
\label{toycase}       
\end{figure}

 The penalization is responsible for part of this noise, as it brings non-linearities to the stiffness, as can be seen on the toy case in Fig. \ref{toycase}. On this figure, the stiffness-volume Pareto fronts for a simple traction problem are represented for two different grids and two different penalizations. The designs corresponding to three points are pictured underneath. Evaluating the same designs with a (post-treatment) penalization $p=1$ instead of the original one $p=3$ gives a smoother Pareto front, and closer to the one that would be obtained with an infinitely fine mesh. Therefore, in order to smooth out our Pareto fronts, we decide to evaluate the compliances/stiffnesses of each design with a penalization of 1 at the end of the optimization process. We however keep the penalization parameter as $p=3$ during optimization, in order to avoid grey designs. 
 
\begin{figure}
\begin{subfigure}{.5\textwidth}
  \includegraphics[width=.99\linewidth]{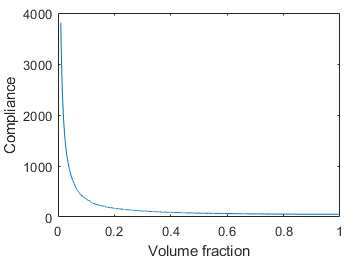}
  \caption{Pareto front}
\end{subfigure}%
\begin{subfigure}{.5\textwidth}
  \includegraphics[width=.99\linewidth]{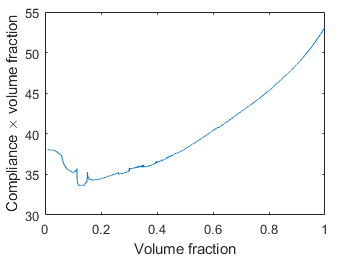}
  \caption{Pareto front multiplied by the volume fraction}
\end{subfigure}%
\caption{Raw compliance - volume fraction Pareto front and Pareto front normalized by the volume fraction.}
\label{mono}
\end{figure}

 When considering the Pareto front obtained with this different final evaluation, we obtain Fig. \ref{mono}a. This seems relatively smooth, but this is not the case in Fig. \ref{mono}b, representing the same Pareto front multiplied by the density. This noisy Pareto front means that different local optimums are attained by different neighboring points. In particular, the fact that the Pareto front is not a monotonous function shows that the global optimum is not reached. This can also be a problem if the Pareto front is to be used in a larger optimization as in Section \ref{MDO}. In this case non-monotony can lead the optimizer to be stuck in local optimums. Therefore, more work is needed to improve the Pareto front obtained.
 
\begin{figure}
\begin{subfigure}{.5\textwidth}
  \includegraphics[width=.99\linewidth]{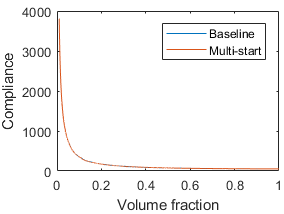}
  \caption{Whole Pareto front}
\end{subfigure}%
\begin{subfigure}{.5\textwidth}
  \includegraphics[width=.99\linewidth]{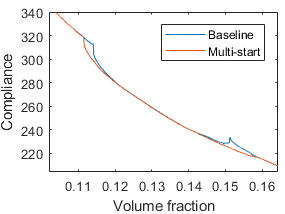}
  \caption{Zoom on noisy zone.}
\end{subfigure}%
\caption{Compliance-volume Pareto fronts obtained with the baseline approach and the multi-start approach.}
\label{multi}
\end{figure}

In order to obtain better optimums, a first idea is to use a multi-start strategy. We use eleven initial designs, depicted in Appendix \ref{initdes}. For each point (i.e. value of volume fraction), we carry out an optimization starting with each of these initial designs. All these initial designs lead to a final design the compliance/stiffness of which is evaluated. We keep the best of these final designs as a point for our Pareto front. The Pareto front obtained is compared to the original Pareto front, named \emph{baseline}, in Fig. \ref{multi}. The initial designs are the same for any problem considered but vary with the volume fraction constraint, as explained in Appendix \ref{initdes}.

\begin{figure}
  \centering
  \includegraphics[width=0.7\textwidth]{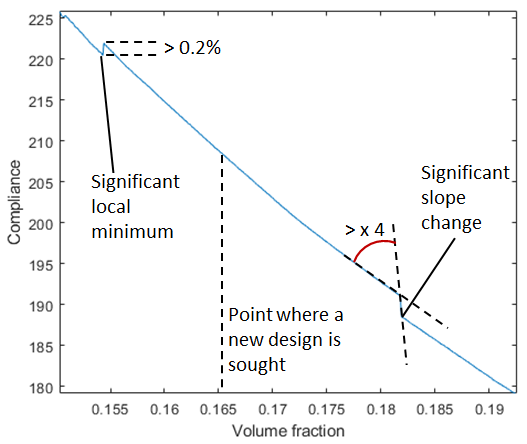}
\caption{Zoom on part of a Pareto front used to define initial points for a better Pareto front.}
\label{signif}       
\end{figure}

There is still some noise. In order to obtain even better Pareto fronts, we decide to go further. Based on the Pareto front obtained with the multi-start strategy, we select some significant local minima, in the compliance times density graph. We define significant local minima, as local minima for which the compliance is lower by at least 0.2\% compared to the compliance of the next point with higher volume fraction. We also select some significant compliance drops in the Pareto front slope. Those two significant points appear in Fig. \ref{signif}. We keep the final designs corresponding to these points. We then build a second Pareto front for which each point is obtained from two initial designs: the closest point corresponding to a local minimum to the left of the point considered (lower volume fraction), and the closest point corresponding to a compliance drop to the right of the point considered (higher volume fraction).

\begin{figure}
  \centering
  \includegraphics[width=0.9\textwidth]{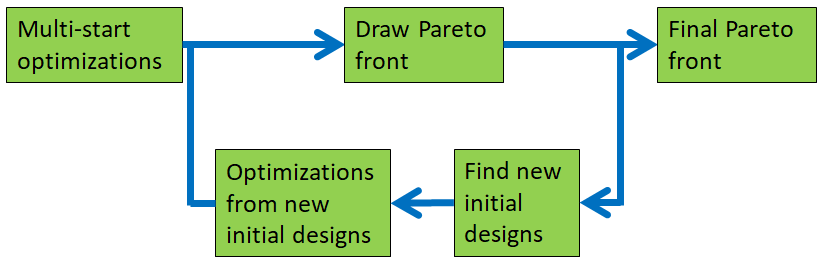}
\caption{Process to improve the Pareto fronts by multiple loops based on increasingly optimal initial designs.}
\label{process}       
\end{figure}

This second Pareto front contains better values than the first, but still exhibits local minimums and drops. The process can therefore be repeated to obtain a third Pareto front and so on. The total process is illustrated in Fig \ref{process}.

\begin{figure}
\begin{subfigure}{.5\textwidth}
  \includegraphics[width=.99\linewidth]{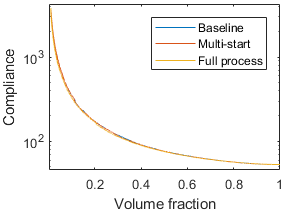}
  \caption{Whole Pareto front}
\end{subfigure}%
\begin{subfigure}{.5\textwidth}
  \includegraphics[width=.99\linewidth]{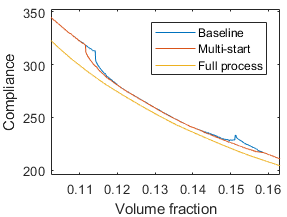}
  \caption{Zoom on noisy zone.}
\end{subfigure}%
\caption{Compliance-volume Pareto fronts obtained with the baseline approach, the multi-start approach or with the full process.}
\label{full}
\end{figure}

This process is computationally intensive, but enables much better Pareto fronts to be obtained, as can be seen in Fig \ref{full} comparing the Pareto front at the end of the process to the one obtained through a single initial design or through the multi-start strategy. 

\begin{figure}
\begin{subfigure}{.5\textwidth}
  \includegraphics[width=.99\linewidth]{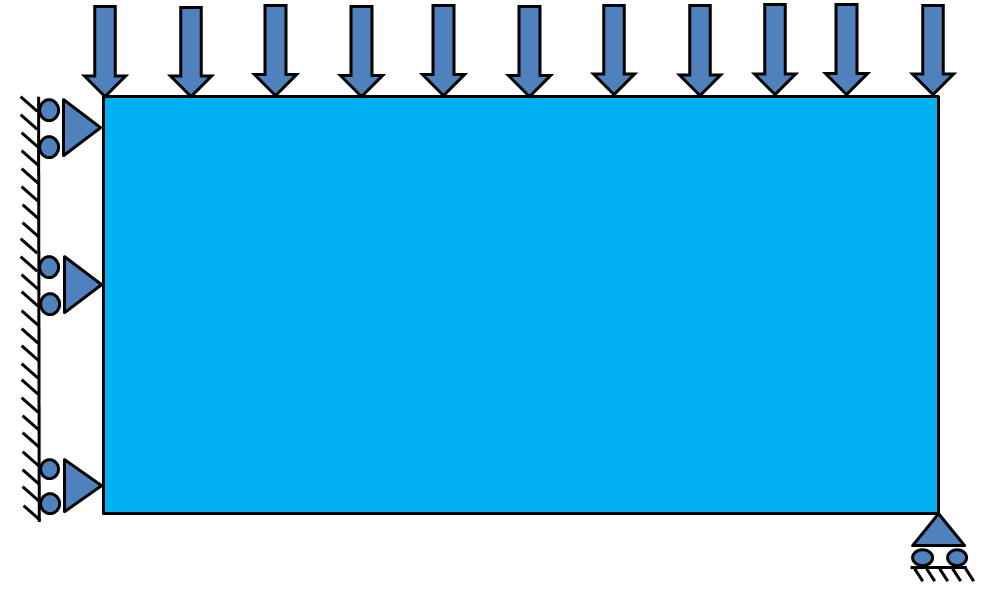}
  \caption{A bridge problem}
\end{subfigure}%
\begin{subfigure}{.5\textwidth}
  \includegraphics[width=.99\linewidth]{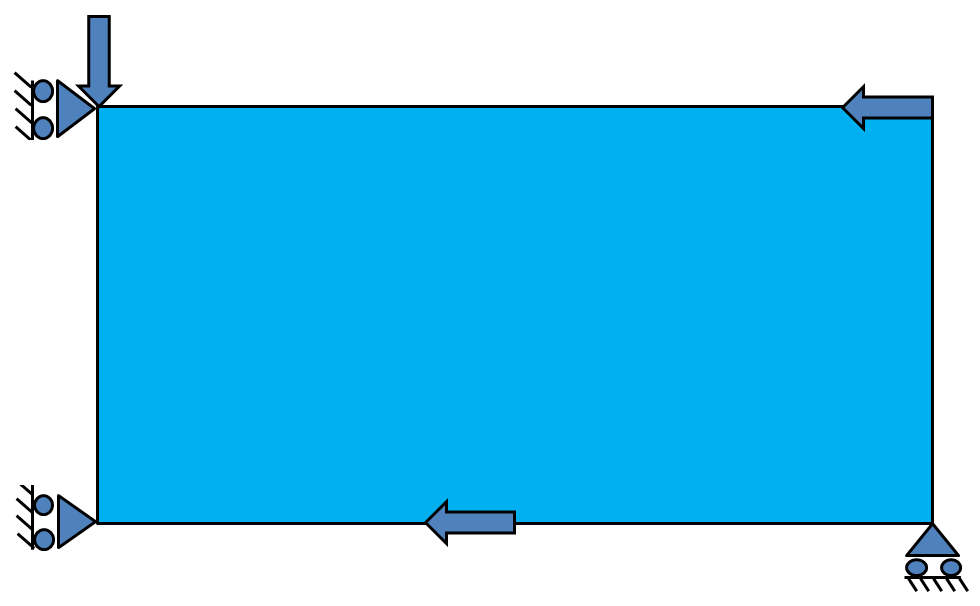}
  \caption{A more complex problem}
\end{subfigure}%
\caption{Two other problems considered to test the versatility of the method.}
\label{probs}
\end{figure}

\begin{figure}
\begin{subfigure}{.5\textwidth}
  \includegraphics[width=.99\linewidth]{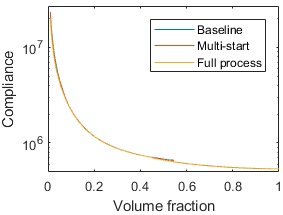}
  \caption{Whole Pareto front}
\end{subfigure}%
\begin{subfigure}{.5\textwidth}
  \includegraphics[width=.99\linewidth]{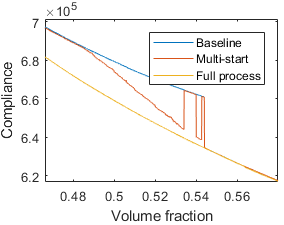}
  \caption{Zoom on noisy zone.}
\end{subfigure}%
\caption{Compliance-volume Pareto fronts obtained with the baseline approach, the multi-start approach or with the full process, for the bridge problem.}
\label{fullbridge}
\end{figure}

\begin{figure}
  \centering
  \includegraphics[width=0.8\textwidth]{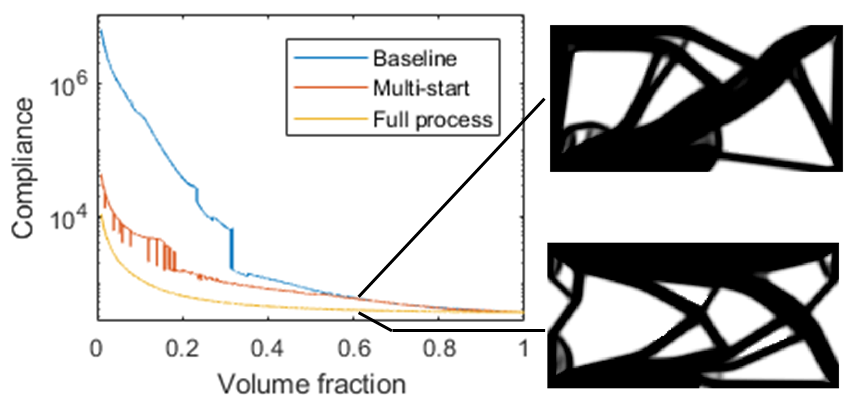}
\caption{Compliance-volume Pareto fronts obtained with the baseline approach, the multi-start approach or with the full process, for the more complex problem. Two designs of same volume fraction are shown.}
\label{fullcomplex}       
\end{figure}

Pareto fronts are also built in the same way for very different problems to test the versatility of this approach. A bridge,  (Fig \ref{probs}a) and a more complex problem (Fig \ref{probs}b) are considered. The respective Pareto fronts appear in Fig \ref{fullbridge} and \ref{fullcomplex}.

\begin{figure}
\begin{subfigure}{.5\textwidth}
  \includegraphics[width=.99\linewidth]{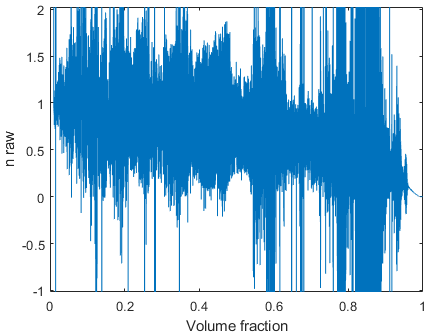}
  \caption{}
\end{subfigure}%
\begin{subfigure}{.5\textwidth}
  \includegraphics[width=.99\linewidth]{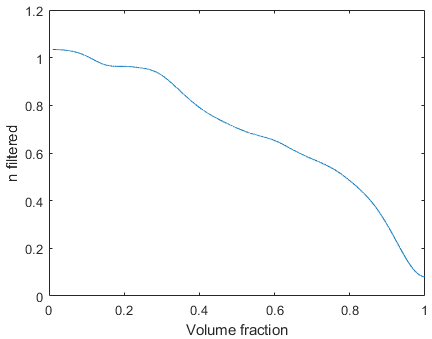}
  \caption{}
\end{subfigure}%
\caption{Raw (a) and filtered (b) ER.}
\label{nfunc}
\end{figure}

With these smoother Pareto fronts, the ER studied in Section \ref{proof} can now be plotted. This is done in Fig. \ref{nfunc}a for the MBB beam. This plot is extremely noisy, as it is based on the derivative of the compliance-volume fraction Pareto front (Eq. (\ref{n1})). This plot is therefore filtered to lead to the plot in Fig. \ref{nfunc}b. This filter is done in two steps. First, we impose the compliance-volume Pareto front to be a decreasing function of $V_f$ by keeping for each point the minimum of all the compliance values with a volume fraction below or equal to the point considered. In the second step, a weighted average of neighbouring points with a Gaussian kernel is applied. A standard deviation of 0.04 is used. The result is much smoother.

\begin{figure}
\begin{subfigure}{.5\textwidth}
  \includegraphics[width=.99\linewidth]{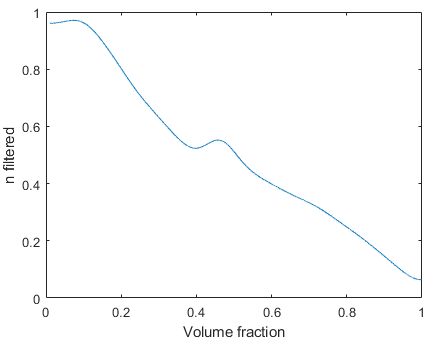}
  \caption{Bridge problem.}
\end{subfigure}%
\begin{subfigure}{.5\textwidth}
  \includegraphics[width=.99\linewidth]{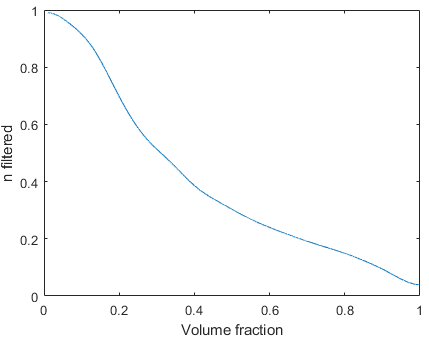}
  \caption{Complex problem.}
\end{subfigure}%
\caption{Filtered ER for the bridge and complex problem.}
\label{nprobs}
\end{figure}

We apply the same filtering approach to the two other problems considered. The smooth ER obtained can be seen in Fig. \ref{nprobs}.

These smooth plots of the ER enable us to check the properties found in Section \ref{proof}. We can indeed see that this function lies in [0,1], tends towards a value of 1 in $V_f=0^+$, and is mostly decreasing (apart from a "bump" for the bridge problem). For the three problems considered, the value of this function tends towards 0 in $V_f=1$ (it appears a bit higher in the filtered plots because of the filtering). This means that the last material added does not bring significant compliance improvement and can be considered useless from this point of view.

\subsection*{A meta-model based on the analysis of the ER}
\label{buildmeta}

In order to carry out material and design simultaneous selection without an excessively heavy computational burden, a simple function representing the Pareto front (i.e. a meta-model) is needed, as explained in Section \ref{MDO}.

We use the properties of the ER obtained in Section \ref{proof} to build a meta-model as simple as possible presenting the same properties. The expression $f_{meta}$ in Eq. (\ref{metamodel}a) is considered to represent the compliance-volume Pareto front $f$. It results in the expression in Eq. (\ref{metamodel}b) for an approximation of the ER.
\begin{subequations}
\begin{equation}
f_{meta}(x)\stackrel{\rm def}{=}A(\frac{1}{x}+Bx^C) 
\end{equation}
Here, A, B and C are tunable real constants and $B>0$
\begin{equation}
n_{meta}(x)=\frac{-xf_{meta}'(x)}{f_{meta}(x)}=\frac{1-BCx^{C+1}}{1+Bx^{C+1}}
\end{equation}
\label{metamodel}
\end{subequations}

\begin{figure}
\begin{subfigure}{.5\textwidth}
  \includegraphics[width=.99\linewidth]{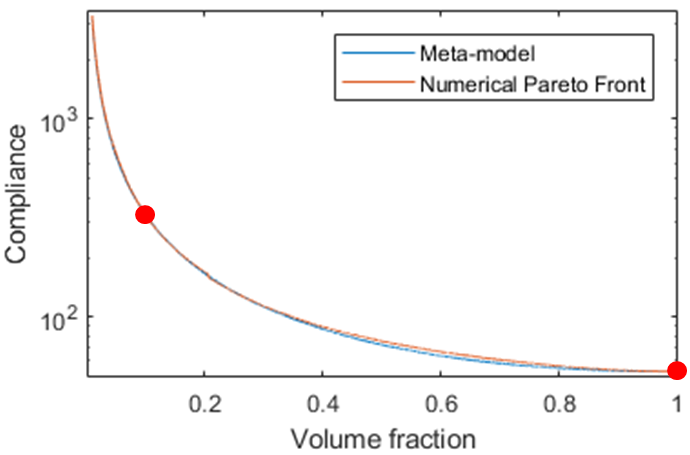}
  \caption{MBB problem}
\end{subfigure}%
\begin{subfigure}{.5\textwidth}
  \includegraphics[width=.99\linewidth]{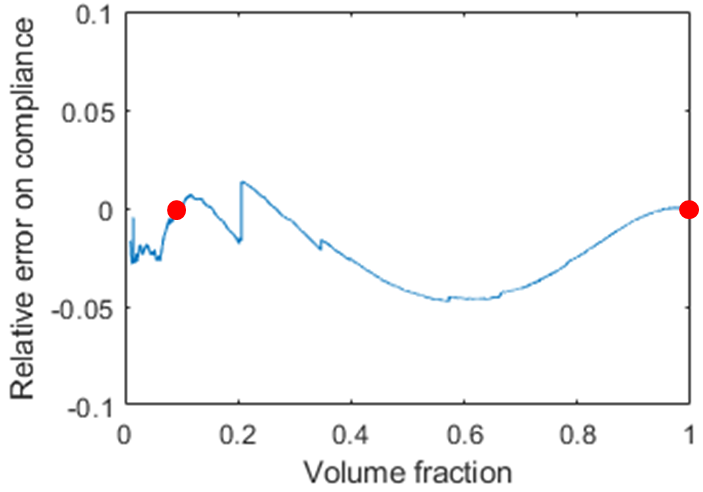}
\end{subfigure}%
\newline
\begin{subfigure}{.5\textwidth}
  \includegraphics[width=.99\linewidth]{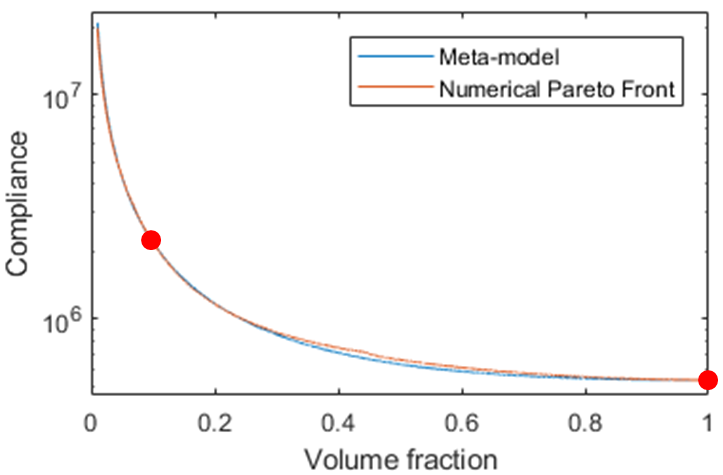}
  \caption{Bridge problem}
\end{subfigure}%
\begin{subfigure}{.5\textwidth}
  \includegraphics[width=.99\linewidth]{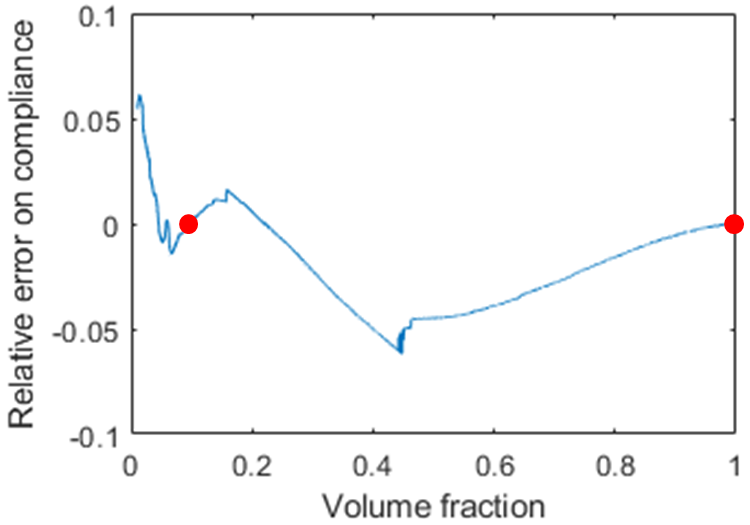}
\end{subfigure}%
\newline
\begin{subfigure}{.5\textwidth}
  \includegraphics[width=.99\linewidth]{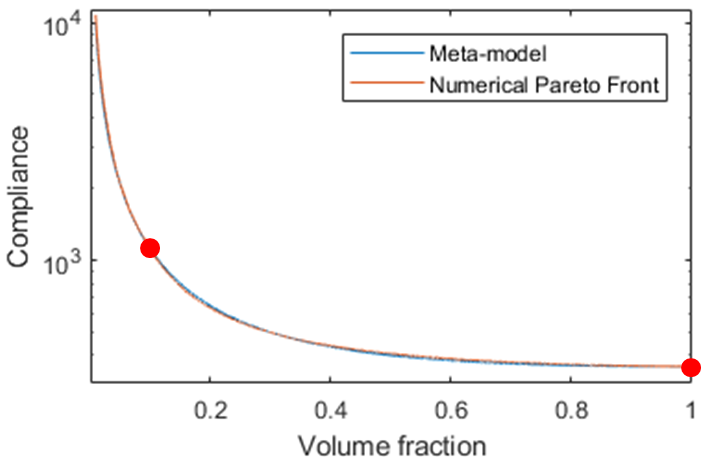}
  \caption{More complex problem}
\end{subfigure}%
\begin{subfigure}{.5\textwidth}
  \includegraphics[width=.99\linewidth]{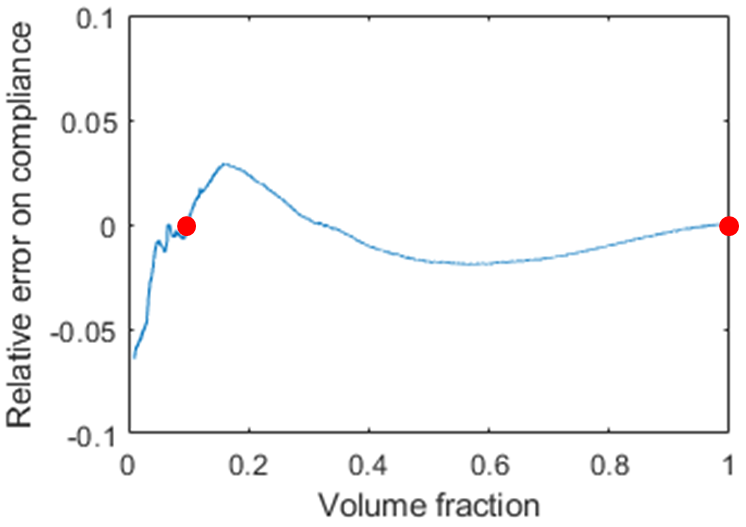}
\end{subfigure}%
\caption{Meta-models superimposed on the numerical Pareto fronts obtained by the full process (a, c, e). The corresponding relative errors are also shown (b, d, f). The red points indicate the points used to fit the meta-model. Only the point at 0.1 is obtained through topology optimization.}
\label{Paretometa}
\end{figure}

By choosing $C=\frac{1}{B}$, we can make sure $n(0)=1$ and $n(1)=0$. In this expression two parameters are left (A and B). As the idea is to obtain the Pareto front meta-model in a computationally efficient way, we can compute only two points to fit the meta-model. The design with a volume fraction of 1 can be obtained directly without needing a topology optimization. It is therefore chosen as one of the two points. The other point is chosen at a volume fraction of 0.1 in order to gather a large range between the two points. Therefore, only one topology optimization is needed to obtain the meta-model of the whole Pareto front. This is done for the three previous examples, and the meta-models obtained are superimposed on the Pareto fronts obtained in Section \ref{numPareto}. We obtain Fig. \ref{Paretometa}. The relative errors are also shown to better illustrate the precision of this meta-model.

The maximum obtained error is 4.7\% for the MBB problem, 6.2\% for the bridge problem and 6.4\% for the more complex problem. This is a very low error, considering that this Pareto front meta-model is appropriate for any compliance-based topology optimization problem, and is obtained through only one topology optimization (at a volume fraction of 0.1). This is very useful for material and design simultaneous selection, discussed in the next section.

\section*{Material and design simultaneous selection for minimum mass}
\label{MDO}

The simultaneous selection of material and 2D design for mass minimization as in Eq. (\ref{galpb}) is discussed in \cite{duriez_ecodesign_2022} for the simple case where the thickness is a free variable. The properties of the ER and the Pareto front meta-model studied in the previous sections enable us to tackle the more complex case where the thickness is imposed.

Ashby's index in this case is $\rho f^{-1}(\frac{U_{max} E t}{F})$ (Eq. (\ref{AshbyIndex})). The material minimizing this index is the optimal material. If we had access to the reciprocal of the compliance-volume fraction Pareto front $f^{-1}$, we could directly compute Ashby's index for every material and select the material with the lowest index. However, $f^{-1}$ is a very complex function, and we do not have access to it unless we compute the whole Pareto front as in Section \ref{numPareto}, which is very computationaly expensive. Therefore, other strategies are described below.

\subsection*{Restricting material choice using the ER}
\begin{figure}
\begin{subfigure}{.5\textwidth}
  \includegraphics[width=.99\linewidth]{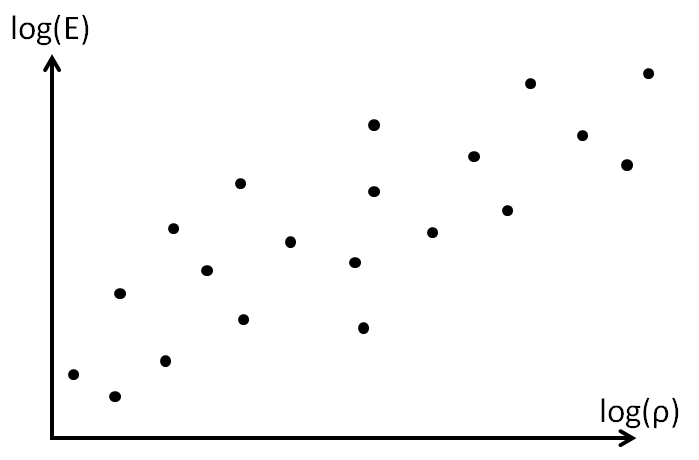}
  \caption{All the materials}
\end{subfigure}%
\begin{subfigure}{.5\textwidth}
  \includegraphics[width=.99\linewidth]{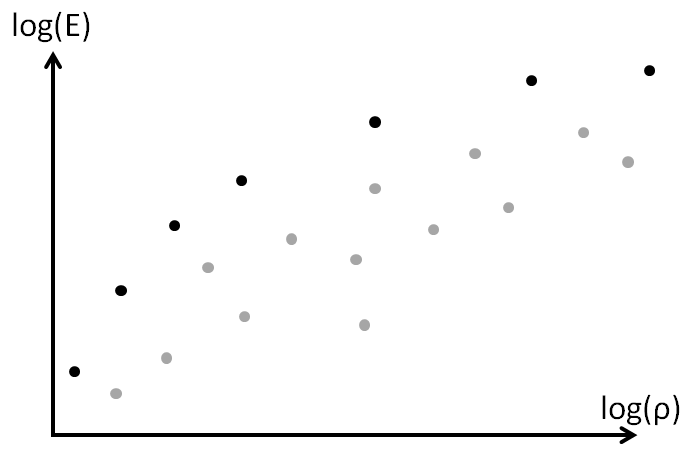}
  \caption{Materials on the Pareto front selected}
\end{subfigure}%
\caption{First screening on a set of fictional materials. Each dot on the Ashby chart corresponds to a specific material.}
\label{ashby}
\end{figure}

There is an enormous number of materials that can be chosen to make the design. A first step is to reduce this number to only materials that could possibly be optimal. As we are looking to minimize mass with a stiffness (displacement) constraint, we want a material with the lowest possible density $\rho$ and the highest possible Young's modulus E. Therefore, only materials on the density-Young's modulus Pareto front can be optimal. We only keep those as candidates. This first screening is illustrated in Fig. \ref{ashby}, representing the Ashby chart (\cite{ashby_materials_2004}) of a set of fictional materials.

To restrict this choice further, we can use the properties of the compliance-volume fraction Pareto front found in Section \ref{proof}. These properties will enable us to prove the following lemma and proposition.

\begin{lem}
\begin{equation}
    \forall a > 1, \forall x \in [f(1), +\infty[, \quad  a f^{-1}(a x) \leq f^{-1}(x)
\end{equation}
with $f^{-1}$ the reciprocal function of the compliance-volume fraction Pareto front $f$.
\label{lem:f-1}
\end{lem}

\begin{proof}

Let g be a function defined by $g(x)=\ln(f^{-1}(x)), \forall x \in [f(1), +\infty[$     (with $f^{-1}(f(1))=1$). Then, its derivative $g^{\prime}$ is such that $g^{\prime}(x)=\frac{{{f^{-1}}^{\prime}}(x)}{f^{-1}(x)}$ with ${f^{-1}}^{\prime}$ as the derivative of $f^{-1}$. The following equations are then verified successively.

\begin{subequations}
\begin{align}
n(x) =
\frac{- x f'(x)}{f(x)} \leq 1 \\
\frac{- f^{-1}(y) f'(f^{-1}(y))}{y} \leq 1 \quad with \quad y=f(x)\\
\frac{f^{-1}(y)}{y {f^{-1}}'(y)} \geq -1 \quad ({f^{-1}}'(y)=\frac{1}{f'(f^{-1}(y))})\\
g'(y) \leq \frac{-1}{y}\\
\int_z^{az}g'(y)dy \leq \int_z^{az} \frac{-1}{y}dy\\
\ln(f^{-1}(az))-\ln(f^{-1}(z)) \leq -\ln(az) + \ln(z) \\
\frac{f^{-1}(ax)}{f^{-1}(x)} \leq \frac{1}{a}
\end{align}
\label{proofafax}
\end{subequations}
\end{proof}

\begin{prop}
 Let us consider the material, defined by its density $\rho_1$ and its Young's modulus $E_1$, with the lowest ratio $\displaystyle{\frac{\rho}{E}}$. Then, the optimal material for the problem stated in Eq. (\ref{mass2}) has a density $\rho_2$ higher than $\rho_1$.
\end{prop}

\begin{proof}
We proceed again by contradiction.
Let us suppose that the optimal material, defined by its density $\rho_2$ and its Young's modulus $E_2$, is such that $\rho_2 < \rho_1$ with $\rho_1$ the density of the material with the lowest $\displaystyle{\frac{\rho_1}{E_1}}$ ratio. As the optimal material is the one with the lowest Ashby's index, 
\begin{equation}
\rho_2 f^{-1}(\frac{U_{max} E_2 t}{F}) < \rho_1 f^{-1}(\frac{U_{max} E_1 t}{F})
\end{equation}
$f^{-1}$ is decreasing, since it is the reciprocal function of a decreasing function $f$. We also know that $\displaystyle{\frac{\rho_1}{E_1}<\frac{\rho_2}{E_2}}$ by the definition of material 1. Therefore,
\begin{equation}
\rho_2 f^{-1}(\frac{U_{max} E_2 t}{F}) < \rho_1 f^{-1}(\frac{U_{max} E_1 t}{F}) < \rho_1 f^{-1}(\frac{U_{max} E_2 t \rho_1}{F \rho_2})
\end{equation}
$\frac{\rho_1}{\rho_2}>1$, therefore, using Lemma \ref{lem:f-1} we have
\begin{equation}
f^{-1}(\frac{U_{max} E_2 t}{F})< \frac{\rho_1}{\rho_2} f^{-1}(\frac{U_{max} E_2 t}{F }\frac{ \rho_1}{\rho_2}) \leq f^{-1}(\frac{U_{max} E_2 t}{F})
\end{equation}

This is an evident contradiction and proves the proposition.
\end{proof}

\begin{figure}
\begin{subfigure}{.48\textwidth}
  \includegraphics[width=.99\linewidth]{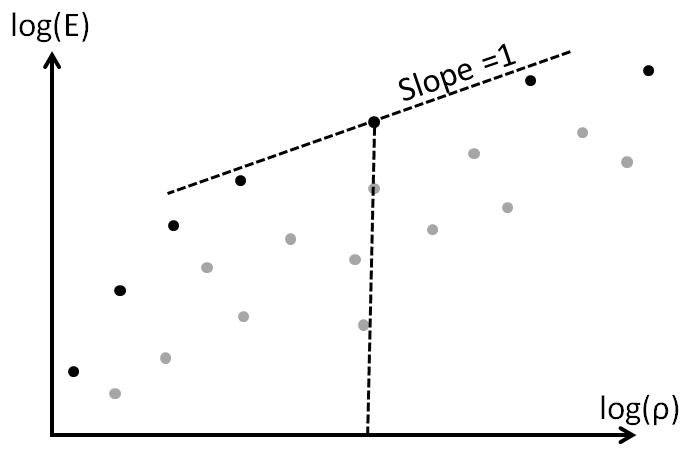}
  \caption{Materials initialy selected. }
\end{subfigure}%
\begin{subfigure}{.48\textwidth}
  \includegraphics[width=.99\linewidth]{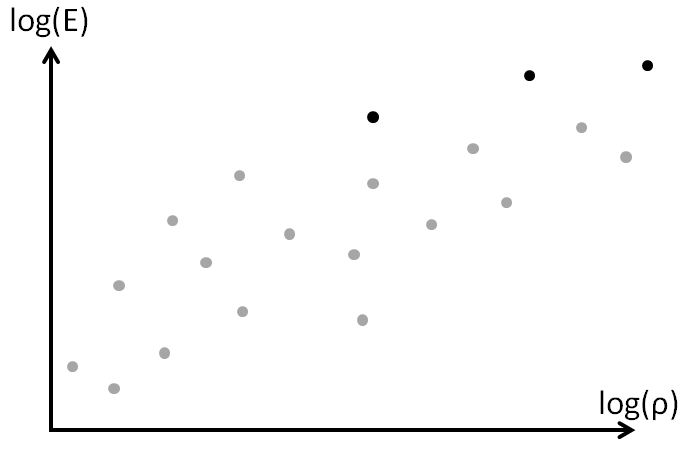}
  \caption{Materials with a density higher than the one with an optimal $\displaystyle{\frac{\rho}{E}}$ ratio.}
\end{subfigure}%
\caption{Second screening on the same set of fictional materials as Fig. \ref{ashby}. The material with the optimal $\displaystyle{\frac{\rho}{E}}$ ratio is located where the line of slope 1 is tangent to the Pareto front in the log-log Ashby chart.}
\label{ashby2}
\end{figure}

We can now use this proposition to restrict the number of optimal material candidates. We can compute the ratio $\displaystyle{\frac{\rho}{E}}$ for every material, and only keep the materials with a higher density than the one with the lowest ratio. It can be counter-intuitive to keep only materials with a high density when we want to minimize mass. In fact, materials with higher density but on the Young's modulus-density Pareto front, have a higher Young's modulus. This means that a lower volume fraction can be used, resulting in a higher ER value: less material will be used more efficiently. This second screening is illustrated in Fig. \ref{ashby2}.

This can be interpreted by looking directly at the ER. As its value stays lower or equal to 1, only materials with the lowest $\displaystyle{\frac{\rho^n}{E}}$ ratio with $n<1$ can be optimal. On the Ashby graph of Fig. \ref{ashby2}, these are the materials on the Pareto front to the left of material 1.

\subsection*{Final choice using the compliance-volume Pareto front meta-model}
\label{choicestrategy}
After the screening steps describe previously, much fewer materials are left. If there is only 1 material left, then it is optimum and the optimal design will be the one along the compliance-volume fraction Pareto front enabling to comply with the stiffness constraint with the lowest volume fraction possible. A first guess for this optimal volume fraction can be found using the reciprocal function of the meta-model of the compliance-volume fraction Pareto front. This first guess should be close to the optimal volume fraction. The compliance can then be computed at this point. If more precision is wanted, a new meta-model can then be made based on this point (according to Eq (\ref{metamodel}) in Section \ref{buildmeta}).
Through this process, a very small number of topology optimizations are necessary to find the optimal volume fraction.

If several optimal material candidates are left after screening, the compliance-volume fraction Pareto front meta-model can be used to find its numerical reciprocal. This approximate reciprocal function is then used to evaluate an approximate Ashby's index of all the candidates left. The material with the lowest approximate Ashby's index is selected and the optimal volume fraction and design can be found as described above. Because these are only approximate Ashby's index, if the best material's index is too close to the second best, then the optimal design and volume fraction can be calculated for these two (or more) candidates in order to choose the best one.

Therefore, the compliance-volume fraction Pareto front meta-model enables a much faster selection of materials, using only a few topology optimizations. Indeed, the alternative would consist in finding the optimal volume fraction for every material on the density-Young's modulus Pareto front and then selecting the one resulting in the lowest mass. This would need several topology optimizations for every material.

\subsection*{Example}
The method described here is illustrated on the same case as in \cite{duriez_ecodesign_2022} for comparison purposes.
The mass of an MBB beam of length 2000mm and height 500mm is to be minimized. A force of 20kN is applied, and a maximum deflection of 5mm is accepted. In the original problem, the thickness was a free variable, leading to an easy solution. Here, we will consider this thickness to be imposed and equal to 5mm. The same four metallic materials are candidates and listed in Tab. \ref{mats}.

\begin{table}[h]
\caption{The different materials considered}
\begin{tabular*}{\hsize}{@{\extracolsep{\fill}}llll@{}}
\toprule
Material & E (GPa) & $\rho$  $(kg/m^3)$ & $\frac{\rho}{E} (\frac{kg}{Nm} .10^{-9})$ \\
\midrule
Aluminum alloy (7475) &  70.8 &  2795 & 39.5\\
Stainless steel (AISI 347)  & 197 & 7915 &  40.2\\
Tinanium alloy (Ti-6Al-4V) &  116 & 4400 &   37.9 \\
Inconel 713 &  205 &  7900 & 38.5 \\
\bottomrule
\end{tabular*}
\label{mats}
\end{table}

We get rid of stainless steel through the first screening. Indeed, it is not on the Pareto front, as Inconel has both a higher Young's modulus and a lower density. We get rid of the aluminum alloy through the second screening. Indeed, it has a lower density than the titanium alloy which has the lowest $\frac{\rho}{E}$ ratio. After screening, we are therefore left with inconel and the titanium alloy.

We can now use the compliance-volume fraction Pareto front meta-model for the MBB beam to find the approximate Ashby indices of these two materials ($f_4$ in Eq. (\ref{AshbyIndex})). These indices appear on Tab. \ref{index}. These indices can be interpreted as the density of a material that would fill the whole design space and result in a part of the same mass.
\begin{table}[h]
\caption{The Ashby indices of the different materials considered, based on the meta-model}
\begin{tabular*}{\hsize}{@{\extracolsep{\fill}}lllll@{}}
\toprule
Material & $E (GPa)$ & $\rho (kg/m^3)$ & $\frac{\rho}{E} (\frac{kg}{Nm} .10^{-9})$ & $f_4 (kg/m^3)$ \\
\midrule
Tinanium alloy (Ti-6Al-4V) &  116 & 4400 &   37.9 & \textcolor{green}{100.1} \\
Inconel 713 &  205 &  7900 & 38.5 & 101.4 \\
\bottomrule
\end{tabular*}
\label{index}
\end{table}

Therefore, Titanium is the best material according to the approximate Ashby index. The resulting volume fraction ($V_f$) is 0.023 and the resulting mass of the part is 0.5kg. Only single topology optimization has been used for this selection process. Iterating on the volume fraction found as in Section \ref{choicestrategy} would give results with a better precision.
If the real Pareto front was used instead of the meta-model, a mass of 0.51kg is found for this material, and it is still optimum. However, obtaining the real Pareto front is much more computationally expensive. 

This very small volume fraction and thickness means the design displays very slender components, and buckling would need to be considered. The very small difference between the indices of the two materials mean that Inconel is almost as good as Titanium in this case.

If the load is multiplied by 20, the optimal material becomes Inconel, based on the real Pareto front. This is explained by the fact that Titanium being lighter and less stiff, the design with this material will have a higher volume fraction. Therefore, material is less efficiently used ($n(V_f)$ is a decreasing function, see Section \ref{proof}). This change of optimal material choice is accurately captured by the meta-model.

\section*{Conclusion}
In this paper, an ``efficiency ratio" was introduced to study compliance-volume fraction Pareto fronts. It was shown that adding more material is less efficient at higher volume fractions. This was used to restrain the number of optimal material candidates. The properties of this ER also helped build a generic meta-model that could easily be fitted to any compliance-volume fraction problem with only one topology optimization, with an acceptable precision (6.4\% error at most). This very efficient meta-model enables fast selection of the optimal material for mass minimization under stiffness constraints on an example. However, in this theoretical work, no limits were applied to the design's members' thicknesses. This can lead to manufacturability or buckling issues. These should be considered in future works.

\subsection*{Data availability}
The datasets generated during and/or analysed during the current study are available from the corresponding author on reasonable request.

\bibliography{sample}

\section*{Acknowledgements}

The authors thank the AMX program for funding this research.

\section*{Author contributions statement}

E.D. and J.M. conceived the material selection method. E.D. and M.C. built the proofs. All authors reviewed the manuscript.

\section*{Competing interests}
The authors declare no competing interests.

\end{document}